\documentclass[final,1p,times]{elsarticle}
	\usepackage{amssymb}
	\usepackage{amsmath}
	\usepackage{graphicx}	
	\usepackage{tikz,pgfplots}
	\pgfplotsset{compat=newest}
	\usepackage{circuitikz}
	\usepackage{pgf}	

\journal{Electric Power System Research}

\begin{document}

\begin{frontmatter}
\title{A Quadratic Approximation for the Optimal Power Flow in Power Distribution Systems}
\author{Alejandro Garces}%
\address{email: alejandro.garces@utp.edu.co \\ Universidad Tecnologica de Pereira. \\ 
AA: 97 - Post Code: 660003 - Pereira, Colombia}

\begin{abstract}

This paper presents a quadratic approximation for the optimal power flow in power distributions systems. The proposed approach is based on a linearized load flow which is valid for power distribution systems including three-phase unbalanced operation. The main feature of the methodology is its simplicity. The accuracy of the proposed approximation is compared to the non-linear/non-convex formulation of the optimal power flow using different optimization solvers. The studies indicate the proposed approximation provides a very accurate solution for systems with a good voltage profile. Results over a set of 1000 randomly generated test power distribution systems demonstrate this solution can be considered for practical purposes in most of the cases. An analytical solution for the unconstrained problem is also developed. This solution can be used as an initialization point for a more precise formulation of the problem.

\end{abstract}
		

\end{frontmatter}

\section{Introduction}

Optimal power flow  (OPF) is a classic problem for transmission system operation which has been extensively studied in the scientific literature \cite{historyopf,Carpentier1979}. The increasing penetration of renewable energies and the possibilities offered by communications in future smart-grids allow the use of OPF in power distribution systems \cite{activereactivedistribution,dynamicopf} and especially in micro-grids \cite{josefguerrero,convexopf1}. 

OPF is a challenging problem due to the high number of nonconvex constraints. Newton-Raphson, descendent gradient and interior points methods are traditionally employed to obtain suboptimal solutions \cite{historyopf}. These methods allow decoupled formulations in the context of transmission networks, since nodal voltages are usually close to $1\angle 0$ and reactance/resistance ratio of transmission lines is frequently high. A good quality initial solution as well as a simple modeling which allows fast calculation of derivatives, are key features for a fast and accurate solution of the problem \cite{efficientimplementation}. 

In power distribution systems, it is a further challenging problem due to the unbalanced operation and low reactance/resistance ratio of distribution lines.  Therefore, heuristic algorithms based on artificial intelligence have been proposed to find good solutions to the problem \cite{carriersopf,biogeographical}.  Evolutionary algorithms \cite{evolutionary,evolutionary3} and particle swarm optimization \cite{swarm} have demonstrated to be efficient approaches for the problem.  This type of algorithms allow an accurate modeling of the system by including constrains otherwise very difficult to consider.  However,  heuristic algorithms do not guarantee optimality and can be computationally cumbersome for real time operation. 

Another approach for the problem in both, transmission and distribution networks, is the use of relaxations and simplifications in order to "convexify" the problem \cite{convexbranch,convexrelaxation,suficientconditions,Equivalentrelaxations,exactconvexrelaxation}.  Semidefinite programming is one of the most promising modeling techniques for this propose \cite{reducedcomplexity,localsolutions}. The main advantage to reformulate a problem as a convex optimization problem is the capability to find global optimal solutions in an efficient way \cite{libroconvexoptimization}.   In addition, a convex formulation allows in some cases, the use of distributed methods.  This is a key feature for future smart-grids.

The difficulty of the OPF lies in the non-convex nature of the load flow equations rather than in the number of variables.  Different convex approximation have been proposed in the literature to address this problem. For example, in \cite{marti} a curve-fitting technique was used in order to linearize voltage-dependent load models. Other analytical approaches were presented in \cite{mitlineal} and \cite{yoflujocarga}. 

This paper introduces a quadratic convex approximation for the OPF in power distribution systems.  This approximation is based on the linear formulation of the power flow presented in \cite{yoflujocarga}.  Different consideration are made ending at a non-iterative analytical solution for the relaxed problem.  Both, the quadratic convex model and the analytical relaxed model are extensible to three-phase unbalanced distribution systems.  These results have many potential applications including:

\begin{itemize}
	\item As initial point for other non-linear or heuristic algorithms
	\item As a practical solution in systems in which a close-to-the-optimal solution is acceptable
	\item In markets regulation where a convex formulation is desired (i.e the solution is unique and do not depend on the used algorithm).
	\item In real time operation were a fast solution is required
	\item As part of other algorithm that requires to call many times an OPF as a sub-routine. 
	\item As sensitivity analysis for power distribution systems
\end{itemize}

Unlike conventional formulations, the proposed approximation uses complex voltages as state variables represented in rectangular form. Although similar formulations have been proposed before \cite{augmentedpf,currentinjection,quadraticopf,evolutionaryquadratic} the results presented here are different in three main aspects:  First, the proposed formulation seeks an approximated model rather than an efficient implementation of a conventional Newton-based algorithm. This approximation has theoretical and practical applications from the power engineering stand point.  Second, modeling and linearization is made entirely on complex variable before split in real and imaginary part for the optimization process.  Off course, it might be possible first split and then linearize, but modeling in complex variables allows a straightforward extension to the three-phase case and inclusion of complex constrains. Third, an non-iterative solution is found for the relaxed case.  Due to the non-convex characteristic of the problem, even the unconstrained case is difficult and can lead to local optimums \cite{localsolutions}. Therefore, a global non-iterative solution is useful even for initialization purposes \cite{resistivenet}.

The remainder of this paper is organized as follows.  Section II reviews the non-linear non-convex formulation of the OPF and analyzes the advantages of a rectangular formulation in power distribution systems.  Section III  presents the quadratic convex approximation of the OPF.  In Section IV, further approximation are considered to end at a non-iterative analytical solution for the relaxed problem.  In Section IV, the methodology is extended to three-phase unbalanced systems.  Finally, Section V presents simulation results performed over an extensive set of test systems before Section VI concludes.

\section{Formulation of the OPF for power distribution systems}

Different formulations for the OPF have been proposed in the scientific literature as a result of contributions from many researchers in this area. Two main formulations can be considered namely Polar-OPF and Rectangular-OPF, according to the representation of the state variables. Both formulations are equivalent. In the first case, decision variables are active and reactive power of distributed generators while voltages are state variables represented in polar form.  In the second case, decision variables are currents injected by generators and state variables are voltages (both represented in rectangular form).  Different objective functions can be considered including minimal generation costs, maximum market surplus and minimum losses, among others. In this paper, the minimum losses OPF is considered although the methodology can be extended for other objective functions.

Rectangular-OPF is less common in the literature than Polar-OPF \cite{historyopf}.  However, it has some advantages in power distribution systems, especially in those cases where distributed generators are operated at constant power factor.  In this formulation, voltages and currents are represented in rectangular form ($(v_{r},v_{i})(i_{r},i_{i})$) as given in Eqs. (\ref{eq:modelorectangular1}) to (\ref{eq:modelorectangular2})

\begin{equation}
\text{Minimize}  \; P_{L} = \left( 2 \sum\limits_{k=1}^{N} g_{(k0)}\cdot v_{r(k)}\cdot v_{(0)}\right) +
  \left(\sum\limits_{k=1}^{N}\sum\limits_{m=1}^{N} g_{(km)} \cdot v_{r(k)} \cdot v_{r(m)}\right)  + \left(\sum\limits_{k=1}^{N}\sum\limits_{m=1}^{N} g_{(km)}\cdot v_{i(k)}\cdot v_{i(m)}\right) 
	\label{eq:modelorectangular1}
\end{equation}

subject to

\begin{equation}
  i_{r(k)} = \:\sum\limits_{m=1}^{N} g_{(km)}\cdot v_{r(m)} - b_{(km)}\cdot v_{i(m)} 
	\label{eq:ld1}
\end{equation}	

\begin{equation}
  i_{i(k)} = \sum\limits_{m=1}^{N} g_{(km)}\cdot v_{i(m)} + b_{(km)}\cdot v_{r(m)} 
\end{equation}

\begin{equation}	
  v_{(k)}^{2}\cdot i_{r(k)} = \left( p_{(k)}\cdot v_{(k)}^{\alpha_{(k)}} + P_{g(k)}\right)\cdot v_{r(k)} + 
               \left( q_{(k)}\cdot v_{(k)}^{\alpha_{(k)}} + Q_{g(k)}\right)\cdot v_{i(k)} 
\end{equation}

\begin{equation}							
  v_{(k)}^{2}\cdot i_{i(k)} = \left(p_{(k)}\cdot v_{(k)}^{\alpha_{(k)}} + P_{g(k)}\right)\cdot v_{i(k)} - 
               \left(q_{(k)}\cdot v_{(k)}^{\alpha_{k}} + Q_{g(k)}\right)\cdot v_{r(k)} 
\label{eq:ld2}							
\end{equation}

\begin{equation}							
 v^{2}_{(k)} 	= v^{2}_{r(k)} + v^{2}_{i(k)} 
\end{equation}

\begin{equation}
  v_{min} \leq v_{(k)}  \leq v_{max} 
\end{equation}

\begin{equation}	
  P_{min} \leq P_{g(k)} \leq P_{max} 
\end{equation}

\begin{equation}	
  Q_{min} \leq Q_{g(k)} \leq Q_{max}   
	\label{eq:modelorectangular2}
\end{equation}

where subscripts $r,i$ represent real and imaginary part, and subscripts $(k),(m)$ represent nodes (with $(0)$ the slack node).  Moreover, $g_{(km)},b_{(km)}$ represent real and imaginary components of the nodal admittance matrix respectively, $p$ and $q$ are the nodal active and reactive power, $P_{g},Q_{g}$ are the active and reactive power delivered by distributed generators, and finally $\alpha$ is an exponent that represent the ZIP model of each load (i.g $0$ for constant power loads, $1$ for constant current and $2$ for constant impedance).  

The main advantage of the Rectangular-OPF is that coupling between nodes are represented by a linear equation while the non-linear/non-convex equations are isolated to each bus.  The main source of non-linear equations are constant power loads. Distributed generators can be considered as PQ buses for mathematical optimization modeling even if operated at constant voltage. Set point of the voltage can be calculated after the optimization is performed. On the other hand, most of the loads in distribution systems require a model which considers constant power, constant current and constant impedance loads (i.e the ZIP model).  Therefore, a linearization of the constant power loads is required in order to obtain a convex approximation.

\section{Quadratic approximation}

This section presents the quadratic approximation from the Rectangular-OPF.  The key step in this development is the linearization of the load flow equations which was first presented in \cite{yoflujocarga}.  For the sake of completeness it is briefly presented below.

Let us consider a power distribution system whose topology is described by the nodal admittance matrix as follows:

\begin{equation}
\left(\begin{array}{c}   \mathbb{I}_{0} \\ \mathbb{I}_{N} \end{array} \right) =
\left(\begin{array}{cc}  \mathbb{Y}_{00} & \mathbb{Y}_{0N} \\ \mathbb{Y}_{N0} & \mathbb{Y}_{NN} \end{array}\right) \cdot
\left(\begin{array}{c}   \mathbb{V}_{0} \\ \mathbb{V}_{N} \end{array}\right)
\end{equation}

where $0$ is the substation node (slack) and $N=\left\{1,2,\dots n\right\}$ are the remained nodes. Along this section, a blackboard bold variable represents a complex matrix or vector while an unbold variable with sub index $r$ or $i$ represents its real or imaginary part respectively.  Slack node is assumed at $0^{o}$ hence $\mathbb{V}_{0}=V_{0}$

Current and voltage in each node $k \in N$ are related to the ZIP model of the loads by Eq. (\ref{eq:corriente})

\begin{equation}
	\mathbb{I}_{k} = \frac{\mathbb{S}^{*}_{Pk}}{\mathbb{V}^{*}_{k}} + \mathbb{S}^{*}_{Ik} + \mathbb{S}^{*}_{Zk}\cdot \mathbb{V}_{k} \label{eq:corriente}
\end{equation}

Notice that component P of the ZIP model is the only non-linear equation.  An approximation of this term was proposed in \cite{yoflujocarga}, based on the Laurent series expansion \cite{librocomplejos} of the function $\mathbb{F}(\mathbb{V})=1/\mathbb{V}$ within a closed region $\left\|  1-\mathbb{V} \right\| < \alpha < 1 $ as follows:  

\begin{equation}
	\frac{\mathbb{S}^{*}_{Pk}}{\mathbb{V}^{*}_{k}} \approx \mathbb{S}^{*}_{Pk}\cdot\left( 2 - \mathbb{V}^{*}_{k}\right)  
\end{equation}

This approximation can be also obtained by spliting Eq. (\ref{eq:corriente}) in real and imaginary part, and expanding in a Taylor series around $\mathbb{V}=1$. However, it is straightforward to maintain the formulation on the complex set. It allows a linear and compact equation for the load flow  which is analog to the DC power flow for power transmission networks as follows:

\begin{equation}
	\mathbb{A} + \mathbb{B} \cdot \mathbb{V}^{*}_{N} + \mathbb{C} \cdot \mathbb{V}_{N} = 0 
	\label{eq:fcds}
\end{equation}

with

\begin{equation}
	\mathbb{A} = \mathbb{Y}_{N0} V_{0} - 2\cdot \mathbb{S}^{*}_{PN} - \mathbb{S}^{*}_{IN} 
\end{equation}
\begin{equation}	
	\mathbb{B} = diag(\mathbb{S}^{*}_{PN}) 
\end{equation}
\begin{equation}	
	\mathbb{C} = \mathbb{Y}_{NN}  - diag(\mathbb{S}^{*}_{ZN})
\end{equation}

It is important to notice that not supposition is made related to whether or not the system is radial or balanced. Unlike the DC formulation of the power flow in transmission networks, nodal voltages are not fixed to 1~pu but obtained by the power flow itself. Thus, PV nodes are not considered. Nevertheless, PV nodes, if any, are less common in power distribution systems than in transmission networks. In addition, distributed generation can be considered as PQ nodes during the optimization process regardless of the type of control in real time operation.

The exactitude of the linear formulation of the power flow and its extension to three-phase systems was discussed in \cite{yoflujocarga} as function of the maximum voltage drop ($\delta_{max}$). 

\begin{equation}
	\left\| 1 - \mathbb{V} \right\| \leq \delta_{max}
	\label{eq:voltagedrop}
\end{equation}

The proposed approximation is valid for power distribution systems in which Eq. (\ref{eq:voltagedrop}) is fulfilled. For example, a $\delta_{max} = 0.3$ (i.e $0.7 < \left\| \mathbb{V} \right\| < 1.3$) results in a maximum expected error of $10\%$ which decreases as $\left\| \mathbb{V} \right\| \rightarrow 1$.   In addition, angle for the slack node must be $0^{o}$.  Otherwise, a rotation is required as presented in Section \ref{sec:trifasico} for three-phase systems.

Let us include the distributed generators in Eq. (\ref{eq:fcds}) as follows:

\begin{equation}
	\mathbb{A} + \mathbb{B} \cdot \mathbb{V}^{*}_{N} + \mathbb{C} \cdot \mathbb{V}_{N} + D \cdot \mathbb{S}^{*}_{G} = 0 \label{eq:linpf}
\end{equation}

where $D=\left( d_{k,m} \right)_{N\times G}$ is a real matrix such that $d_{km} = 1$ if the generator $m$ is connected to the node $k$, and 0 otherwise.  This affine equation replaces the non-linear expressions given  by Eqs. \ref{eq:ld1} to \ref{eq:ld2}.  The additional convex constrains required to complete the optimization model are given in Eq.
(\ref{eq:fcoq}):

\begin{equation}
\begin{array}{rc}
		\text{Minimize} & P_{L} = Real\left( \mathbb{V}^{T}_{N} \cdot \mathbb{Y}_{NN} \mathbb{V}^{*}_{N} + 2\cdot \mathbb{V}_{N}^{T}      \mathbb{Y}_{N0}\cdot \mathbb{V}_{0} \right) \\ \\
 \text{subject to} & 

		\mathbb{A} + \mathbb{B} \cdot \mathbb{V}^{*}_{N} + \mathbb{C} \cdot \mathbb{V}_{N} + D \cdot \mathbb{S}^{*}_{G} = 0 \\
&   \left\| 1-\mathbb{V} \right\| \leq \delta_{max}  \\
&   \left\| \mathbb{S}_{G} \right\| \leq S_{G(max)} \\
\end{array}
\label{eq:fcoq}
\end{equation}

After separating in real and imaginary part, all equality constrains are affine functions whereas all inequalities are convex function. Symbol $\left\| . \right\|$ could represent absolute-value norm or Euclidean norm of the complex number.  Both cases result in a convex optimization problem.  However, the first case is easier to solve since it is a quadratic/convex optimization problem with linear constrains. 

The optimization model given by Eq. (\ref{eq:fcoq}) is compact and clear due to the representation in complex variables.  This is important for implementation purposes in a software that allows complex variables (i.g Matlab/Octave). Objective function cannot be complex on account of it is not an ordered field.  However, feasible region and decision variables can be efficiently represented on the complex field. Nevertheless, a more conventional representation of this model is obtained by separating in real and imaginary part starting from Eq. (\ref{eq:linpf}) as follows:

\begin{equation}
	\left( \begin{array}{c} D \cdot S_{r} \\ -D \cdot S_{i} \end{array} \right)
	+
	\left( \begin{array}{c} A_{r} \\  A_{i} \end{array} \right)  
	+ 	
	\left( \begin{array}{cc} B_{r} + C_{r}  &  B_{i}-C_{i} \\ 
	                         B_{i} + C_{i}  & -B_{r}+C_{r}  
  \end{array} \right)	
	\cdot 
	\left( \begin{array}{r} V_{r} \\ V_{i} \end{array} \right)
	= 0
	\label{eq:fcee}
\end{equation}

Notice this expression allows to obtain $V_{r}$ and $V_{i}$ as function of $S_{r}$, $S_{i}$ (ie. the active and reactive power delivered by distributed generators) since remaining terms are all constant.  Let us define a constant matrix $M_{N\times N}$ as follows:

\begin{equation}
	M = \left( \begin{array}{cc} B_{r} + C_{r}  & B_{i}-C_{i} \\ B_{i} + C_{i} & -B_{r}+C_{r}  \end{array} \right)
\end{equation}

then (\ref{eq:fcee}) can be re-written as

\begin{equation}
	\left( \begin{array}{c}
		V_{r} \\ V_{i}
	\end{array}\right) = -M^{-1} \cdot	
	\left( \begin{array}{c}
		A_{r} \\ A_{i}
	\end{array}\right)
	-M^{-1}
	\left(
	\begin{array}{c}
		D\cdot S_{r} \\ -D\cdot S_{i}
	\end{array}\right)	
\end{equation}

let us define $\mathbb{U}=U_{r} + j\cdot U_{i}$ as the vector voltages of the system without distributed generation (i.e $\mathbb{S}_{G}=0$).

\begin{equation}
\left( \begin{array}{c} U_{r} \\ U_{i} \end{array} \right) = 
M^{-1} \cdot
\left( \begin{array}{c} -A_{r} \\  -A_{i} \end{array} \right)
\end{equation}

and $W_{2N\times 2N}$ a real matrix given by

\begin{equation}
	\left( \begin{array}{rr} W_{rr} & W_{ri} \\ W_{ir} & W_{ii}  \end{array} \right) =
	M^{-1} \cdot
	\left( \begin{array}{cc} -D  & 0 \\ 0 & D  \end{array} \right) \label{eq:ws}
\end{equation}

These definitions allow to represent nodal voltages as function of the generated power.

\begin{equation}
	\left( \begin{array}{c} V_{r} \\ V_{i} \end{array} \right) = 
	\left( \begin{array}{c} U_{r} \\ U_{i} \end{array} \right) + 
	\left( \begin{array}{cc} W_{rr} & W_{ri} \\ W_{ir} & W_{ii} \end{array} \right)\cdot
	\left( \begin{array}{c} S_{r} \\ S_{i} \end{array} \right)
	\label{eq:frms}
\end{equation}

Nodal voltages $\mathbb{V}$ are equal to $\mathbb{U}$ when generated power $\mathbb{S}_{G}$ is zero, which correspond with the definition of $\mathbb{U}$. 
At this point Eq. (\ref{eq:frms}) is just another representation of Eq. (\ref{eq:linpf})that is replaced on the optimization model (\ref{eq:fcoq}) as follows:

\begin{equation}
\begin{array}{rc}
		\text{Minimize} &  V^{T}_{r} \cdot {G}_{N} \cdot V_{r} + 
		                          V^{T}_{i} \cdot {G}_{N} \cdot V_{i} + 2 \cdot V^{T}_{r} \cdot G_{0}\cdot V_{0} \\				
 \text{subject to} & \\

		& V_{r} = U_{r} + W_{rr}\cdot S_{r} + W_{ri}\cdot S_{i} \\
		& V_{i} = U_{i} + W_{ir}\cdot S_{r} + W_{ii}\cdot S_{i} \\
		& \left(1-V_{r(k)}\right)^{2} + V_{i(k)}^{2} \leq \delta_{max}^{2} \\
		& S_{r(k)} \leq S_{r(max)} \\
		& S_{i(k)} \leq S_{i(max)} \\
\end{array}
\label{eq:fcoqnuevo}
\end{equation}

where

\begin{eqnarray}
G_{N} = Real(\mathbb{Y}_{NN}) \\
G_{0} = Real(\mathbb{Y}_{N0})
\end{eqnarray}

The model is convex and hence, it can be solved easily no mater the size of the system. $G_{N}$ is a square matrix but $W$is not. Therefore, $V_{r}=V_{i}=0$ are not feasible solutions.  In addition, if maximum drop constraint (Eq. \ref{eq:voltagedrop}) is relaxed, then the model becomes in a simple quadratic optimization with linear constraints which can be solved by well known toolboxes such as quadprog in matlab. Nevertheless, Eq \ref{eq:voltagedrop} must be checked after the optimization since the exactitude of the load flow depends highly of the this constraint. 

\subsection{Analytical relaxed model}

Inequality constraints in Eq. \ref{eq:fcoqnuevo} can be relaxed in order to obtain an analytical relaxed model.  This model can be used for initialization of the non-linear optimal power flow as will be presented in Section \ref{sec:results}.  

Let us include load flow equations into the objective function as follows:

\begin{equation}
\begin{split}
	\text{Minimize} \; \left(U_{r} + W_{rr}\cdot S_{r} + W_{ri}\cdot S_{i} \right)^{T}\cdot G_{N}\cdot \left(U_{r} + W_{rr}\cdot S_{r} + W_{ri}\cdot S_{i} \right)  + \\
	\left(U_{i} + W_{ir}\cdot S_{r} + W_{ii}\cdot S_{i}\right)^{T} \cdot G_{N} \cdot 
	\left(U_{i} + W_{ir}\cdot S_{r} + W_{ii}\cdot S_{i}\right) + \\
	2 \left(U_{r} + W_{rr}\cdot S_{r} + W_{ri}\cdot S_{i} \right)^{T} \cdot G_{0}\cdot V_{0}
\end{split}	
\end{equation}

which can be re-written as

\begin{equation}
	\text{Minimize} \; \frac{1}{2} S^{T} \cdot H \cdot S + F^{T} \cdot S \\
\end{equation}

with

\begin{equation}
	H = \left(\begin{array}{c|c}
  H_{rr} & H_{ri} \\ \hline H_{ir} & H_{ii}
	\end{array} \right)
\end{equation}

\begin{eqnarray}
	  H_{rr} = W^{T}_{rr}\cdot G_{N} \cdot W^{T}_{rr} + W^{T}_{ir}\cdot G_{N} \cdot W^{T}_{ir} \\ 
		H_{ri} = W^{T}_{rr}\cdot G_{N} \cdot W^{T}_{ri} + W^{T}_{ir}\cdot G_{N} \cdot W^{T}_{ii} \\ 
		H_{ir} = W^{T}_{ri}\cdot G_{N} \cdot W^{T}_{rr} + W^{T}_{ii}\cdot G_{N} \cdot W^{T}_{ir} \\
		H_{ii} = W^{T}_{ri}\cdot G_{N} \cdot W^{T}_{ir} + W^{T}_{ii}\cdot G_{N} \cdot W^{T}_{ii}
\end{eqnarray}

and $F = \left( F_{r}, \; F_{i}\right)^{T}$ is given by

\begin{eqnarray}
F_{r}	=  W^{T}_{rr}\cdot G_{N} \cdot U^{T}_{r} + W^{T}_{ir}\cdot G_{N} \cdot U^{T}_{i} + W^{T}_{rr}\cdot G_{0}\cdot V_{0} \label{eq:fr} \\ 
F_{i} =	W^{T}_{ri}\cdot G_{N} \cdot U^{T}_{r} + W^{T}_{ii}\cdot G_{N} \cdot U^{T}_{i} + W^{T}_{ri}\cdot G_{0}\cdot V_{0} 
\label{eq:fi}
\end{eqnarray}

where angle of the slack node is assumed to be equal to zero (i.e $V_{i0}=0$).  The optimal solution of the relaxed optimization problem (i.e analytical relaxed model)  is given by Eq. (\ref{eq:analitico}).

\begin{equation}
	S = -H^{-1} \cdot F \label{eq:analitico}
\end{equation}

This analytical solution has many applications.  In a general case, it can be used as an initial solution for a non-linear optimization methodology or a sensitivity equation for other optimization problems. 
 
\section{Extension to three-phase systems} \label{sec:trifasico}

The proposed approximations can be extended to three-phase unbalanced distribution systems with some additional considerations. Nodal admittance matrix must consider 
a three-phase modeling with a three-phase Slack bus.  As a consequence, $\mathbb{V}_{0}$ is not a constant but a $3\times 1$ vector.  In addition, Y-connected as well as $\Delta$-connected loads must be considered.  Per unit representation is not an advantage in this circumstances. Thus, a new constant $\eta=1/v_{nom}$ is defined in order to permit the linearization around 1.  

The basic structure of Eq. \ref{eq:linpf} remains but now the matrices are redefined as follows:

\begin{equation}
	\mathbb{A} = \mathbb{Y}_{N0}\cdot{V}_{0}-2\eta\cdot J^{T}_{NN}\cdot\mathbb{S}_{PN}^{*}\circ\mathbb{T}_{N}-\eta\cdot J^{T}_{NN}\cdot\mathbb{S}_{IN}^{*}\circ\mathbb{T}_{N}
\end{equation}

\begin{equation}
	\mathbb{B} = \eta^{2}\cdot J^{T}_{NN}\cdot diag(\mathbb{S}_{PN}^{*}\circ\mathbb{T}^{2}_{N})\cdot J_{NN}
\end{equation}

\begin{equation}
	\mathbb{C} = \mathbb{Y}_{NN}-\eta^{2}\cdot J^{T}_{NN}\cdot diag(\mathbb{S}_{ZN}^{*})\cdot J_{NN}
\end{equation}

\begin{equation}
	\mathbb{D} = \eta\cdot J^{T}_{NN} \cdot diag(\mathbb{T}_{N})\cdot D \label{eq:wwe}
\end{equation}

where $\mathbb{T}_{(N\times 1)}$ is a complex vector whose elements correspond to each node such that $\mathbb{T}_{k} = e^{j\phi}$ with $\phi \in \left\{0,-2\pi/3,2\pi/3 \right\}$ according to the sequence. Symbol ($\circ$) is the Hadamard product (i.g the matlab element-wise multiplication) in order differentiate from the conventional matrix multiplication. Moreover,  $J$ is a real matrix which indicates the type of connection of each load: for Y-connected loads, it is an identity matrix while for $\Delta$-connected loads it is a matrix that converts phase voltage into line voltages.  

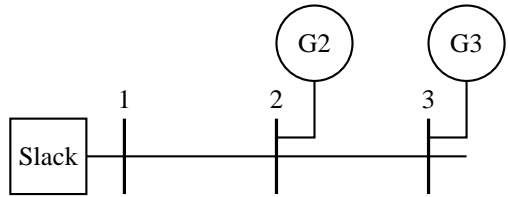
\begin{figure}[htb]
\centering
\begin{tikzpicture}[x=1mm, y = 1mm]

\draw[-,thick] (-5,-5) rectangle +(10,10);
\node at (0,0) {Slack};
\draw[-,thick] (5,0) -- +(50,0);
\draw[-, very thick] (10,-5) -- +(0,10) node[above] {1};
\draw[-, very thick] (30,-5) -- +(0,10) node[above] {2};
\draw[-, very thick] (50,-5) -- +(0,10) node[above] {3};

\draw[-,thick] (35,15) circle (5);
\node at (35,15) {G2};
\draw[-,thick] (35,10) |- (30,2.5);

\draw[-,thick] (55,15) circle (5);
\node at (55,15) {G3};
\draw[-,thick] (55,10) |- (50,2.5);
\end{tikzpicture}
\caption{Three-phase distribution system with distributed generators}
\label{fig:sistematrifasico}
\end{figure}

On the other hand, $D$ has the same definition as in the single-phase case but a new complex matrix $\mathbb{D}$ requires to be defined as Eq. (\ref{eq:wwe}).  Additional constrains related to the three-phase modeling could be included directly on $\mathbb{D}$. for example, some generators could be inappropriate for unbalanced generation (i.g renewable resources integrated by three-legs voltage-source-converters). In those cases, generated power on each phase must be equal despite the unbalanced operation of the rest of the system. This is represented by Eq. (\ref{eq:dsde}).

\begin{equation}
	\mathbb{S}_{k} = \mathbb{S}_{m} = \mathbb{S}_{n}  \label{eq:dsde}
\end{equation}

where $k,m,n$ are different phases of the same three-phase bus.  This constraint can be included directly in D.  To take a simple example, consider a three-phase distribution system with 3 nodes and 2  distributed generators connected to node 2 and 3 as depicted in Fig.\ref{fig:sistematrifasico}.

Generator connected to node 2 can generate three-phase balanced power while generator connected to node 3 can generate power independently to each phase.  The matrix $D$ for this system is given by Eq. (\ref{eq:matd}) where nodes are ordered by phase (i.e fist the phase $a$ for nodes 1,2,3, then phase $b$ and finally phase $c$).

\begin{equation}	
  D = \left(
	\begin{array}{c|ccc}
    0      & 0   & 0 & 0  \\
		1/3    & 0   & 0 & 0  \\
    0      & 1   & 0 & 0 \\
		\hline\\
		0      & 0 & 0 & 0  \\
    1/3    & 0 & 0  & 0  \\
    0      & 0 & 1  & 0  \\
		\hline\\
    0      & 0 & 0 & 0  \\		
    1/3    & 0 & 0 & 0  \\
		0      & 0 & 0 & 1  \\
  \end{array} 
	\right)  \label{eq:matd}
\end{equation}

If generators is and loads are $\Delta$-connected, then $J_{NN}$ is given by

\begin{equation}
	J = \left( \begin{array}{rrr|rrr|rrr}
	1	&  0	&  0	& -1	&  0	&  0	&  0	&  0	&  0  \\
	0	&  1	&  0	&  0	& -1	&  0	&  0	&  0	&  0  \\
	0	&  0	&  1	&  0	&  0	& -1	&  0	&  0	&  0  \\ \hline
	0	&  0	&  0	&  1	&  0	&  0	& -1	&  0	&  0  \\
	0	&  0	&  0	&  0	&  1	&  0	&  0	& -1	&  0  \\
	0	&  0	&  0	&  0	&  0	&  1	&  0	&  0	& -1  \\ \hline
 -1	&  0	&  0	&  0	&  0	&  0	&  1	&  0	&  0  \\
	0	& -1	&  0	&  0	&  0	&  0	&  0	&  1	&  0  \\
	0	&  0	& -1	&  0	&  0	&  0	&  0	&  0	&  1  	
	\end{array}\right)
\end{equation}

As $\mathbb{D}$ is complex, Eqs \ref{eq:ws}, \ref{eq:fr} and \ref{eq:fi} require to be modified as follows

\begin{equation}
	\left( \begin{array}{rr} W_{rr} & W_{ri} \\ W_{ir} & W_{ii}  \end{array} \right) =
	M^{-1} \cdot
	\left( \begin{array}{cc} D_{r}  & D_{i} \\ D_{i} & -D_{r}  \end{array} \right) \label{eq:ws}
\end{equation}

\begin{eqnarray}
\nonumber F_{r}	=  W^{T}_{rr}\cdot G_{N} \cdot U^{T}_{r} + W^{T}_{ir}\cdot G_{N} \cdot U^{T}_{i} + W^{T}_{rr}\cdot G_{0}\cdot V_{r0} \\  +  W^{T}_{ir}\cdot G_{N}\cdot V_{i0}\\ 
\nonumber F_{i} =	W^{T}_{ri}\cdot G_{N} \cdot U^{T}_{r} + W^{T}_{ii}\cdot G_{N} \cdot U^{T}_{i} + W^{T}_{ri}\cdot G_{0}\cdot V_{r0} \\ +  W^{T}_{ii}\cdot G_{N}\cdot V_{i0}
\label{eq:fr}
\end{eqnarray}

Notice that $V_{i0} \in \bold{R}^{3}$ is not zero in this case.  Both, the quadratic convex approximation and the analytical relaxed model can be extended for three-phase modeling. The problem increases in size but its solution is straightforward.

\section{Results} \label{sec:results}

Different simulations were performed in order to test the accuracy of the proposed methodology. 

The IEEE37 power distribution test system \cite{testfeeders} was used to evaluate the performance of the proposed methodology using a balanced equivalent modeling.  Three generators were included in nodes 708, 732 and 744.  Three optimization models were considered: Rectangular-OPF, quadratic programming formulation and relaxed analytical solution. 
Rectangular-OPF was implemented in GAMS \cite{gams} as a non-linear/non-convex problem. The quadratic programming formulation was solved using the optimization toolbox of Matlab while the relaxed analytical formulation was solved by implementing Eq.(\ref{eq:analitico}) directly in Matlab.  Results for $S_{r(max)}=S_{i(max)}=0.9$ are given in Table \ref{tab:resultadosirrestricto}.  Generated power for the quadratic programming formulation were exactly the same as for the relaxed analytical formulation (i.e inequality constrains were never active) on account of the high capacity of distributed generators.  Moreover, this result was very close to the exact non-convex/non-linear programming formulation. Differences between the exact and the approximated model were less than $3\%$.  Total active power loss of the system without distributed generation was $P_{L}=0.0298$ while the optimized system was $P_{L}=0.0089$. A considerable reduction. 

\begin{table}[htb]
\centering
\caption{Comparative results of the OPF on the IEEE37 test system with $P_{max}=Q_{max}=0.9$ }
\label{tab:resultadosirrestricto}
\begin{tabular}{|c|c|c|c|}
\hline
												& Rectangular-OPF  & Quadratic convex & Analytical relaxed \\
		 \textbf{Variable}  &  (GAMS)          & model (MATLAB)   & model (MATLAB)     \\           
\hline \hline
	$\mathbb{S}_{708}$ & 0.849+j0.400 & 0.8580+j0.4153 & 0.8580+j0.4153 \\
	$\mathbb{S}_{732}$ & 0.042+j0.028 & 0.0420+j0.0210 & 0.0420+j0.0210 \\
	$\mathbb{S}_{744}$ & 0.368+j0.191 & 0.3881+j0.1913 & 0.3881+j0.1913 \\
	$P_{L}$            & 0.00876             & 0.0089         & 0.0089 \\
\hline		
\end{tabular}
\end{table}   

A second simulation was performed on the same test system but this time the total capacity of generation was decreased to $P_{max}=Q_{max}=0.7$. Results are given in Table \ref{tab:conrestricciones}.  Relaxed Analytical solution was not as good as in the previous case on account of the activation of inequality constrains.  Nevertheless, the quadratic programming formulation was still a very good approximation with an error of less than $3\%$.  

\begin{table}[htb]
\centering
\caption{Comparative results of the OPF on the IEEE37 test system with $P_{max}=Q_{max}=0.7$}
\label{tab:conrestricciones}
\begin{tabular}{|c|c|c|}
\hline
												& Rectangular-OPF & Quadratic convex  \\
		 \textbf{Variable}  & (GAMS)          & model (MATLAB)     \\           
\hline \hline
	$\mathbb{S}_{708}$ & 0.700+j0.400 & 0.700+j0.415\\
	$\mathbb{S}_{732}$ & 0.147+j0.028 & 0.154+j0.021\\
	$\mathbb{S}_{744}$ & 0.409+j0.190 & 0.413+j0.191\\
	$P_{L}$            & 0.0088        & 0.00896 \\
	
\hline		
\end{tabular}
\end{table}   

The analytical relaxed model as initialization methodology for non-linear programming algorithms was also evaluated on the IEEE37 test system. Table \ref{tab:inicializacion} shows the number of iterations with and without initialization for different solvers available in GAMS. The number of iterations was drastically reduced in some cases (for example in IPOPT solver).

\begin{table}[htb]
\centering
\caption{Comparative results of the OPF on the IEEE37 test system using different solvers }
\label{tab:inicializacion}
\begin{tabular}{|c|c|c|}
\hline
                        & With             & Without         \\ 
		 \textbf{Solver}    & Initialization   & Initialization   \\           
\hline \hline
  MINOS   & 4   & 6  \\
  KNITRO  & 6   & 9  \\
  CONOPT  & 16  & 17 \\
	IPOPT   & 9   & 75 \\	
	LINDO   & 18  & 50 \\		
	PATHNLP & 1   & 47 \\
\hline		
\end{tabular}
\end{table} 

Accuracy of the proposed approximation was analyzed in a more general context by a set of simulations over 1000 randomly generated radial distribution systems with parameters given in Table \ref{tab:random}.  Test system were compelled to be radial, not as a condition of the methodology, but as manner to simply the algorithm for random generation of test cases.

\begin{table}[htb]
\centering
\caption{Parameters for the generation of random distribution systems}
\label{tab:random}
\begin{tabular}{|l|c|c|}
\hline											
		 \textbf{Parameter}  & \textbf{Min} & \textbf{Max}  \\           
\hline \hline
	Number of nodes ($n$) & 30 & 60 \\
	Line resistance & 0.001 & 0.0170 \\
	Line inductive impedance & 0.001 & 0.0170 \\
	Line capacitive admittance & 0.000 & 0.0002 \\
	Load power & 0.000 & 0.2000 \\
	Load power factor & 0.7 & 1.0 \\
	Constant power loads & $50\%$ of $n$ & \\		
	Distributed generators & $5\%$ of $n$ & \\
\hline		
\end{tabular}
\end{table}  

First, Total power loss was calculated on each test system without distributed generation using a back-forward sweep algorithm \cite{Renato}.  Then, an OPF was estimated using the proposed methodology. Next, a new power flow was calculated using the exact model.  After all these calculations were performed on each randomly generated test system, a histograms were drawn to analyze the results.

\begin{figure}[htb]
\footnotesize
\centering
\begin{tikzpicture}
\begin{axis}[scale only axis,width=7.3cm, height=3.5cm,xmajorgrids,ymajorgrids, xlabel={Reduction of power loss $[\%]$}, ylabel ={Frequency}, ymax = 200]
\addplot [ybar, fill=blue!50!green, solid, draw=black, thick, bar width=18pt,nodes near coords,
nodes near coords align={vertical}]
coordinates{
    	(	   4.8076  ,  37.0000	)
	    ( 	 14.3614 ,  44.0000	)
	    (	   23.9151 ,  59.0000	)
	    (	   33.4689 ,  65.0000	)
	    (	   43.0227 , 105.0000	)
	    (	   52.5764 , 125.0000	)
	    (	   62.1302 , 142.0000	)
	    (	   71.6840 , 172.0000	)
	    (	   81.2377 , 172.0000	)
	    (	   90.7915 ,  79.0000	)
};
\end{axis}
\end{tikzpicture}
\caption{Histogram the improvement in terms of power losses for 1000 randomly generated test power distribution systems}
\label{fig:mejora}
\end{figure}
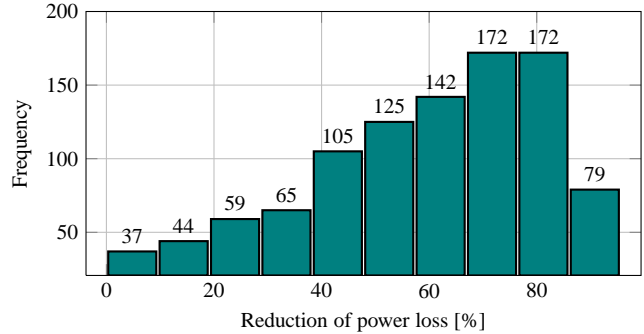

Figure \ref{fig:mejora} depicts the improvement in terms of power losses between the systems without distributed generation and the systems optimized by the proposed approximation of the OPF. In most of the cases, improvement was more than $50\%$.  As expected, distribution generation improves the losses.

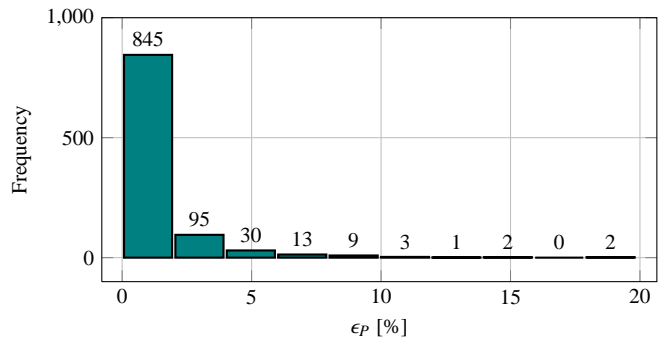
\begin{figure}[htb]
\footnotesize
\centering
\begin{tikzpicture}
\begin{axis}[scale only axis,width=7.3cm, height=3.5cm,xmajorgrids,ymajorgrids, xlabel={$\epsilon_{P}\;[\%]$}, ylabel ={Frequency}, ymax = 1000]
\addplot [ybar, fill=blue!50!green, solid, draw=black, thick, bar width=18pt,nodes near coords,
nodes near coords align={vertical}]
coordinates{
    (	    1.0003 , 845.0000	)
	  (	    2.9854 ,  95.0000	)
    (	    4.9705 ,  30.0000	)
    (	    6.9556 , 13.0000	)
    (	    8.9407 ,   9.0000	)
    (	   10.9258 ,   3.0000	)
    (	   12.9109 ,   1.0000	)
    (	   14.8960 ,   2.0000	)
    (	   16.8810 ,        0	)
    (	   18.8661 ,   2.0000	)
};
\end{axis}
\end{tikzpicture}
\caption{Histogram the error of power losses between the analytical relaxed methodology of OPF and the exact model for 1000 randomly generated test power distribution systems}
\label{fig:perdidas}
\end{figure}

\begin{figure}[htb]
\footnotesize
\centering
\begin{tikzpicture}
\begin{axis}[scale only axis,width=7.3cm, height=3.5cm,xmajorgrids,ymajorgrids, xlabel={$\epsilon_{V}\;[\%]$}, ylabel ={Frequency}, ymax = 1000]
\addplot [ybar, fill=blue!50!green, solid, draw=black, thick, bar width=18pt,nodes near coords,
nodes near coords align={vertical}]
coordinates{
    (	    0.5014 , 838.0000	)
    (	    1.4807 , 109.0000	)
    (	    2.4600 ,  27.0000	)
    (	    3.4394 ,  16.0000	)
    (	    4.4187 ,   4.0000	)
    (	    5.3980 ,   2.0000	)
    (	    6.3773 ,   2.0000	)
    (	    7.3567 ,        0	)
    (	    8.3360 ,   1.0000	)
    (	    9.3153 ,   1.0000	)
};
\end{axis}
\end{tikzpicture}
\caption{Histogram the maximum voltage error between the analytical relaxed methodology of OPF and the exact model for 1000 randomly generated test power distribution systems}
\label{fig:errorvoltages}
\end{figure}
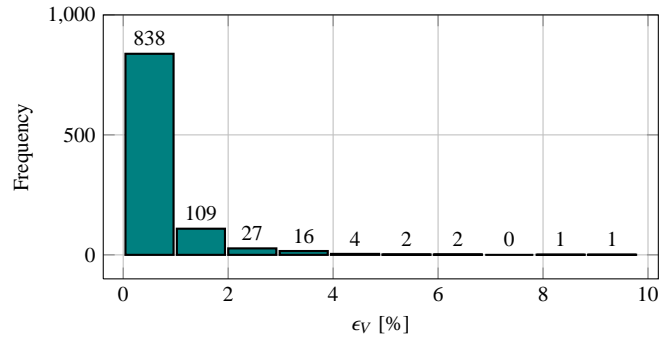

Accuracy of the solutions was demonstrated by an histogram of the error in power losses ($\epsilon_{P}$) and voltages ($\epsilon_{V}$) as shown in Figs \ref{fig:perdidas} and \ref{fig:errorvoltages}.  Power error was less than $5\%$ and voltage error was less than $2\%$ in more than $80\%$ of the cases. Notice that $2\%$ of error does not mean $2\%$ of voltage drop in the system. In these cases, the proposed approximation could be accepted as the solution of the OPF, whereas in the other cases it can be used as good initial approximation for a more accurate optimization model. 

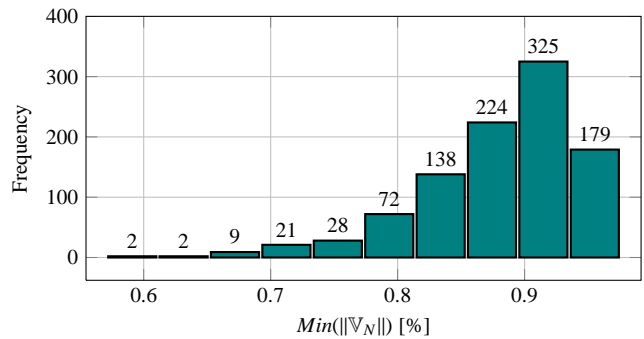
\begin{figure}[htb]
\footnotesize
\centering
\begin{tikzpicture}
\begin{axis}[scale only axis,width=7.3cm, height=3.5cm,xmajorgrids,ymajorgrids, xlabel={$Min(\left\|\mathbb{V}_{N}\right\|)\;[\%]$}, ylabel ={Frequency}, ymax = 400]
\addplot [ybar, fill=blue!50!green, solid, draw=black, thick, bar width=18pt,nodes near coords,
nodes near coords align={vertical}]
coordinates{
   	(	    0.5910  ,  2.0000	)
	  (	    0.6314  ,  2.0000	)
	  (	    0.6719  ,  9.0000	)
	  (	    0.7124  , 21.0000	)
	  (	    0.7529  , 28.0000	)
	  (	    0.7934  , 72.0000	)
	  (	    0.8339  ,138.0000	)
	  (	    0.8743  ,224.0000	)
	  (	    0.9148  ,325.0000	)
	  (	    0.9553  ,179.0000	)

};
\end{axis}
\end{tikzpicture}
\caption{Histogram of the minimal voltage calculated by the back-forward sweep algorithm for 1000 randomly generated test power distribution systems}
\label{fig:minvoltages}
\end{figure}

High errors in other cases can be explained by their loading conditions.  In those cases, the minimal voltage along the feeder is lower than $0.7$ as shown in Fig \ref{fig:minvoltages}.   Fortunately, most of the power distribution systems (as well as the randomly generated test systems) have a good voltage profile and hence the proposed approximation is valid.   

Finally, the three-phase formulation of the optimal power flow was tested on the IEEE37 unbalanced power distribution test system \cite{testfeeders}.  Balanced generation was imposed. Results are summarized in Table \ref{tab:casotrifasico}.   Comparison in terms of maximum voltage error and losses error demonstrate the accuracy of the methodology even in the three-phase case.

\begin{table}[htb]
\centering
\caption{Results of the OPF on the three-phase case}
\label{tab:casotrifasico}
\begin{tabular}{|c|c|}
\hline
  Parameter & value \\
\hline\hline 		
	$\mathbb{S}_{708}$   & $288+137j\; [\text{kW}]$   \\ 
	$\mathbb{S}_{732}$   & $13.8+6.96j\; [\text{kW}]$  \\
	$\mathbb{S}_{744}$   & $131+61j\; [\text{kW}]$  \\
	Losses system base ($\mathbb{S}_{G}=0$) & $0.0335$ \\
	Losses on the optimized system & $0.0122$ \\
	Maximum voltage error $Max(\epsilon{V})$ & $4.7619\times 10^{-4}$ \\
	Losses error $(\epsilon{P})$ & $1.1\%$ \\
\hline	
\end{tabular}
\end{table}
    
\section{Conclusions}

A quadratic approximation of the OPF for power distribution systems was proposed.  This approximation is based on the linear power flow for distribution systems.  Convexity guarantees not only a unique solution but a fast, easy and efficient implementation in commercial optimization toolboxes. Results are very close to the Rectangular-OPF which is a non-linear/non-convex problem.

An analytical solution for the unconstrained problem was also developed. This solution can be used as initial point for the non-linear programming formulation.  Simulations using different commercial solvers available in GAMS shown the importance of a good initial point to accelerate convergence.

Simulation results over a large set of random generated test distribution systems demonstrated the accuracy of the approximation. In most of the cases, the approximated solution had a very low error compared to the back-forward swing algorithm.

The methodology was extended to three-phase systems.  Additional constrains such balanced generation was strait forward included in the methodology.  Y-connected loads as well as $\Delta$-connected loads and single phase generations were included in the model.

The use of a complex expansion instead of a Taylor series for linearization of the system, although equivalent for most of practical purposes, can allow an easy an efficient modeling and implementation of the algorithms.

\bibliographystyle{elsarticle-num}
\bibliography{BibliografiaOPF}

\begin{thebibliography}{10}
\expandafter\ifx\csname url\endcsname\relax
  \def\url#1{\texttt{#1}}\fi
\expandafter\ifx\csname urlprefix\endcsname\relax\def\urlprefix{URL }\fi
\expandafter\ifx\csname href\endcsname\relax
  \def\href#1#2{#2} \def\path#1{#1}\fi

\bibitem{historyopf}
A.~C. Mary B.~Cain, Richard P.~O'Neill, History of Optimal Power Flow and
  Formulations, Federal energy regulatory commision, 2012.

\bibitem{Carpentier1979}
J.~Carpentier, Optimal power flows, International Journal of Electrical Power
  \& Energy Systems 1~(1) (1979) 3 -- 15.

\bibitem{activereactivedistribution}
A.~Gabash, P.~Li, Active-reactive optimal power flow in distribution networks
  with embedded generation and battery storage, Power Systems, IEEE
  Transactions on 27~(4) (2012) 2026--2035.
\newblock \href {http://dx.doi.org/10.1109/TPWRS.2012.2187315}
  {\path{doi:10.1109/TPWRS.2012.2187315}}.

\bibitem{dynamicopf}
S.~Gill, I.~Kockar, G.~Ault, Dynamic optimal power flow for active distribution
  networks, Power Systems, IEEE Transactions on 29~(1) (2014) 121--131.
\newblock \href {http://dx.doi.org/10.1109/TPWRS.2013.2279263}
  {\path{doi:10.1109/TPWRS.2013.2279263}}.

\bibitem{josefguerrero}
Y.~Levron, J.~Guerrero, Y.~Beck, Optimal power flow in microgrids with energy
  storage, Power Systems, IEEE Transactions on 28~(3) (2013) 3226--3234.
\newblock \href {http://dx.doi.org/10.1109/TPWRS.2013.2245925}
  {\path{doi:10.1109/TPWRS.2013.2245925}}.

\bibitem{convexopf1}
E.~Dall'Anese, H.~Zhu, G.~Giannakis, Distributed optimal power flow for smart
  microgrids, Smart Grid, IEEE Transactions on 4~(3) (2013) 1464--1475.
\newblock \href {http://dx.doi.org/10.1109/TSG.2013.2248175}
  {\path{doi:10.1109/TSG.2013.2248175}}.

\bibitem{efficientimplementation}
Q.~Jiang, G.~Geng, C.~Guo, Y.~Cao, An efficient implementation of automatic
  differentiation in interior point optimal power flow, Power Systems, IEEE
  Transactions on 25~(1) (2010) 147--155.
\newblock \href {http://dx.doi.org/10.1109/TPWRS.2009.2030286}
  {\path{doi:10.1109/TPWRS.2009.2030286}}.

\bibitem{carriersopf}
M.~Moeini-Aghtaie, A.~Abbaspour, M.~Fotuhi-Firuzabad, E.~Hajipour, A decomposed
  solution to multiple-energy carriers optimal power flow, Power Systems, IEEE
  Transactions on 29~(2) (2014) 707--716.
\newblock \href {http://dx.doi.org/10.1109/TPWRS.2013.2283259}
  {\path{doi:10.1109/TPWRS.2013.2283259}}.

\bibitem{biogeographical}
A.~Bhattacharya, P.~Chattopadhyay, Application of biogeography-based
  optimisation to solve different optimal power flow problems, Generation,
  Transmission Distribution, IET 5~(1) (2011) 70--80.
\newblock \href {http://dx.doi.org/10.1049/iet-gtd.2010.0237}
  {\path{doi:10.1049/iet-gtd.2010.0237}}.

\bibitem{evolutionary}
E.~Amorim, S.~Hashimoto, F.~Lima, J.~Mantovani, Multi objective evolutionary
  algorithm applied to the optimal power flow problem, Latin America
  Transactions, IEEE (Revista IEEE America Latina) 8~(3) (2010) 236--244.
\newblock \href {http://dx.doi.org/10.1109/TLA.2010.5538398}
  {\path{doi:10.1109/TLA.2010.5538398}}.

\bibitem{evolutionary3}
L.~Honorio, A.~da~Silva, D.~Barbosa, L.~Delboni, Solving optimal power flow
  problems using a probabilistic alpha constrained evolutionary approach,
  Generation, Transmission Distribution, IET 4~(6) (2010) 674--682.
\newblock \href {http://dx.doi.org/10.1049/iet-gtd.2009.0208}
  {\path{doi:10.1049/iet-gtd.2009.0208}}.

\bibitem{swarm}
T.~Niknam, M.~Narimani, J.~Aghaei, R.~Azizipanah-Abarghooee, Improved particle
  swarm optimisation for multi-objective optimal power flow considering the
  cost, loss, emission and voltage stability index, Generation, Transmission
  Distribution, IET 6~(6) (2012) 515--527.
\newblock \href {http://dx.doi.org/10.1049/iet-gtd.2011.0851}
  {\path{doi:10.1049/iet-gtd.2011.0851}}.

\bibitem{convexbranch}
M.~Farivar, S.~Low, Branch flow model: Relaxations and convexification: Part i,
  Power Systems, IEEE Transactions on 28~(3) (2013) 2554--2564.
\newblock \href {http://dx.doi.org/10.1109/TPWRS.2013.2255317}
  {\path{doi:10.1109/TPWRS.2013.2255317}}.

\bibitem{convexrelaxation}
R.~Madani, S.~Sojoudi, J.~Lavaei, Convex relaxation for optimal power flow
  problem: Mesh networks, Power Systems, IEEE Transactions on 30~(1) (2015)
  199--211.
\newblock \href {http://dx.doi.org/10.1109/TPWRS.2014.2322051}
  {\path{doi:10.1109/TPWRS.2014.2322051}}.

\bibitem{suficientconditions}
D.~Molzahn, B.~Lesieutre, C.~DeMarco, A sufficient condition for global
  optimality of solutions to the optimal power flow problem, Power Systems,
  IEEE Transactions on 29~(2) (2014) 978--979.
\newblock \href {http://dx.doi.org/10.1109/TPWRS.2013.2288009}
  {\path{doi:10.1109/TPWRS.2013.2288009}}.

\bibitem{Equivalentrelaxations}
S.~Bose, S.~Low, T.~Teeraratkul, B.~Hassibi, Equivalent relaxations of optimal
  power flow, Automatic Control, IEEE Transactions on PP~(99) (2014) 1--1.
\newblock \href {http://dx.doi.org/10.1109/TAC.2014.2357112}
  {\path{doi:10.1109/TAC.2014.2357112}}.

\bibitem{exactconvexrelaxation}
L.~Gan, N.~Li, U.~Topcu, S.~Low, Exact convex relaxation of optimal power flow
  in radial networks, Automatic Control, IEEE Transactions on 60~(1) (2015)
  72--87.
\newblock \href {http://dx.doi.org/10.1109/TAC.2014.2332712}
  {\path{doi:10.1109/TAC.2014.2332712}}.

\bibitem{reducedcomplexity}
M.~Andersen, A.~Hansson, L.~Vandenberghe, Reduced-complexity semidefinite
  relaxations of optimal power flow problems, Power Systems, IEEE Transactions
  on 29~(4) (2014) 1855--1863.
\newblock \href {http://dx.doi.org/10.1109/TPWRS.2013.2294479}
  {\path{doi:10.1109/TPWRS.2013.2294479}}.

\bibitem{localsolutions}
W.~Bukhsh, A.~Grothey, K.~McKinnon, P.~Trodden, Local solutions of the optimal
  power flow problem, Power Systems, IEEE Transactions on 28~(4) (2013)
  4780--4788.
\newblock \href {http://dx.doi.org/10.1109/TPWRS.2013.2274577}
  {\path{doi:10.1109/TPWRS.2013.2274577}}.

\bibitem{libroconvexoptimization}
S.~Boyd, L.~Vandenberhe, Convex optimization, Cambridge university press, 2004.

\bibitem{marti}
J.~Marti, H.~Ahmadi, L.~Bashualdo, Linear power-flow formulation based on a
  voltage-dependent load model, Power Delivery, IEEE Transactions on 28~(3)
  (2013) 1682--1690.
\newblock \href {http://dx.doi.org/10.1109/TPWRD.2013.2247068}
  {\path{doi:10.1109/TPWRD.2013.2247068}}.

\bibitem{mitlineal}
S.~Bolognani, S.~Zampieri, On the existence and linear approximation of the
  power flow solution in power distribution networks, Power Systems, IEEE
  Transactions on PP~(99) (2015) 1--10.
\newblock \href {http://dx.doi.org/10.1109/TPWRS.2015.2395452}
  {\path{doi:10.1109/TPWRS.2015.2395452}}.

\bibitem{yoflujocarga}
A.~Garces, A linear three-phase load flow for power distribution systems, Power
  Systems, IEEE Transactions on PP~(99) (2015) 1--2.
\newblock \href {http://dx.doi.org/10.1109/TPWRS.2015.2394296}
  {\path{doi:10.1109/TPWRS.2015.2394296}}.

\bibitem{augmentedpf}
A.~Exposito, E.~Ramos, Augmented rectangular load flow model, Power Systems,
  IEEE Transactions on 17~(2) (2002) 271--276.
\newblock \href {http://dx.doi.org/10.1109/TPWRS.2002.1007892}
  {\path{doi:10.1109/TPWRS.2002.1007892}}.

\bibitem{currentinjection}
V.~da~Costa, N.~Martins, J.~Pereira, Developments in the newton raphson power
  flow formulation based on current injections, Power Systems, IEEE
  Transactions on 14~(4) (1999) 1320--1326.
\newblock \href {http://dx.doi.org/10.1109/59.801891}
  {\path{doi:10.1109/59.801891}}.

\bibitem{quadraticopf}
Y.~Tao, A.~Meliopoulos, Optimal power flow via quadratic power flow, in: Power
  Systems Conference and Exposition (PSCE), 2011 IEEE/PES, 2011, pp. 1--8.
\newblock \href {http://dx.doi.org/10.1109/PSCE.2011.5772563}
  {\path{doi:10.1109/PSCE.2011.5772563}}.

\bibitem{evolutionaryquadratic}
S.~Sivasubramani, K.~Swarup, Sequential quadratic programming based
  differential evolution algorithm for optimal power flow problem, Generation,
  Transmission Distribution, IET 5~(11) (2011) 1149--1154.
\newblock \href {http://dx.doi.org/10.1049/iet-gtd.2011.0046}
  {\path{doi:10.1049/iet-gtd.2011.0046}}.

\bibitem{resistivenet}
C.~Tan, D.~Cai, X.~Lou, Resistive network optimal power flow: Uniqueness and
  algorithms, Power Systems, IEEE Transactions on 30~(1) (2015) 263--273.
\newblock \href {http://dx.doi.org/10.1109/TPWRS.2014.2329324}
  {\path{doi:10.1109/TPWRS.2014.2329324}}.

\bibitem{librocomplejos}
F.~Flanigan, Complex Variables, Dover Books on Mathematics, 2010.

\bibitem{testfeeders}
W.~Kersting, Radial distribution test feeders, Power Systems, IEEE Transactions
  on 6~(3) (1991) 975--985.
\newblock \href {http://dx.doi.org/10.1109/59.119237}
  {\path{doi:10.1109/59.119237}}.

\bibitem{gams}
R.~E. Rosenthal, GAMS:General Algebraic Modeling System, GAMS Development
  Corporation, Washington, DC, USA.

\bibitem{Renato}
R.~Cespedes, New method for the analysis of distribution networks, Power
  Delivery, IEEE Transactions on 5~(1) (1990) 391--396.
\newblock \href {http://dx.doi.org/10.1109/61.107303}
  {\path{doi:10.1109/61.107303}}.

\end{thebibliography}

\end{document}